\documentclass[11pt]{amsart}
\usepackage{amsmath,amsfonts,amssymb,amsthm}
\usepackage{latexsym,bm,graphicx}
\usepackage{mathrsfs}
\usepackage{color}
\usepackage{hyperref}

\title[Optimal asymptotic volume ratio for noncompact 3-manifolds]{Optimal asymptotic volume ratio for noncompact 3-manifolds with asymptotically nonnegative Ricci curvature and a uniformly positive scalar curvature lower bound}
\author{Xian-Tao Huang}
\address{School of Mathematics Sun Yat-sen University Guangzhou 510275}
\email{hxiant@mail2.sysu.edu.cn; huangxt28@mail.sysu.edu.cn}
\author{Shuai Liu}
\address{School of Mathematics Sun Yat-sen University Guangzhou 510275}
\email{liush328@mail2.sysu.edu.cn}

\newtheorem{thm}{Theorem}[section]
\newtheorem{prop}[thm]{Proposition}
\newtheorem{lem}[thm]{Lemma}
\newtheorem{cor}[thm]{Corollary}

\newtheorem{conj}[thm]{Conjecture}
\newtheorem{prob}[thm]{Question}

\theoremstyle{definition}
\theoremstyle{remark}

\newtheorem{defn}[thm]{Definition}
\newtheorem{rem}[thm]{Remark}

\numberwithin{equation}{section}

\begin{document}
\maketitle

\begin{abstract} In this paper, we study 3-dimensional complete non-compact Riemannian manifolds with asymptotically nonnegative Ricci curvature and a uniformly positive scalar curvature lower bound.
Our main result is that, if this manifold has $k$ ends and finite first Betti number, then it has at most linear volume growth, and furthermore, if the negative part of Ricci curvature decays sufficiently fast at infinity, then we have an optimal asymptotic volume ratio  $\limsup_{r\rightarrow\infty}\frac{\mathrm{Vol}(B(p, r))}{r}\leq4k\pi$.
In particular, our results apply to 3-dimensional complete non-compact Riemannian manifolds with nonnegative Ricci curvature and a uniformly positive scalar curvature lower bound.


\vspace*{5pt}

\end{abstract}

\section{Introduction}

In this paper, we investigate the influence of positive scalar curvature on the geometric properties of 3-dimensional complete non-compact Riemannian manifolds with asymptotically nonnegative Ricci curvature.
\begin{defn}
We say a Riemannian manifold having asymptotically nonnegative Ricci curvature if for some fixed point $p \in M$, there exists a positive non-increasing continuous function $f(r)$ with $\displaystyle\lim_{r\rightarrow\infty}f(r)=0$, so that for $x\in M$, it holds that $\mathrm{Ric}(x) \geq -f(d(x,p))$.
\end{defn}

It is believed that the uniformly positive lower bound of scalar curvature has an influence on the geometric properties of a manifold.
In \cite{G86}, Gromov raised the following conjecture:
\begin{conj}\label{prob1.1}
Let $(M^n,g)$ be an $n$-dimensional complete non-compact Riemannian manifold with nonnegative Ricci curvature and $\mathrm{Sc}_g\geq1$. Do we have
\begin{equation}\label{1.1}
\limsup_{r\rightarrow\infty}\frac{\mathrm{Vol}(B(p,r))}{r^{n-2}}<\infty?
\end{equation}
\end{conj}

Gromov's conjecture is solved in some special cases.

In \cite{Z22}, Zhu proves that, (\ref{1.1}) holds on $(M^{n}, g)$ with nonnegative Ricci curvature and $\mathrm{Sc}_g\geq1$ provided $M$ has a positive lower bound of the injective radius. 
In the 3-dimensional case, there are more progresses.
Zhu (\cite{Z22}) verified (\ref{1.1}) holds on 3-manifolds with nonnegative Ricci curvature and $\mathrm{Sc}_g\geq1$ under the additional assumption that $\mathrm{inf}_{q\in M}\mathrm{Vol}(B(q,1))\geq v_{0}> 0$.
Later on, Munteanu and Wang (see \cite{MW22}) solved completely the 3-dimensional case by analyzing certain properties of the Green's function. More precisely, the following result has been established in \cite{MW22}:
\begin{thm}[\cite{MW22}]\label{thm1.2}
Let $(M^3,g)$ be a $3$-dimensional complete noncompact Riemannian manifold with nonnegative curvature and $\mathrm{Sc}_g\geq1$.
Then there exists a universal constant $C > 0$ such that for any fixed point $p\in M$ and $r>0$,
\begin{equation}
\frac{\mathrm{Vol}(B(p,r))}{r}<C.
\end{equation}
\end{thm}
We remark that in \cite{MW22} Munteanu and Wang also study $3$-manifolds with certain type of asymptotically nonnegative Ricci curvature.

Later, Chodosh, Li and Stryker (\cite{CLS23}) use the theory of $\mu$-bubbles and Cheeger-Colding's almost splitting theorem to recover Theorem \ref{thm1.2}, with the exception that the constant $C$ in \cite{CLS23} is not universal, but depends on $M$ (on the other hand, we remark that both in \cite{MW22} \cite{CLS23}, the authors consider the growth rate of volume of geodesic balls in the case that the positive scalar curvature lower bound decays compared to some power of the distance function, and the class handled in \cite{CLS23} is wider than that of \cite{MW22}).
In a recent paper \cite{Wang24}, by a more detailed exploration of the method in \cite{CLS23}, Wang is able to recover Theorem \ref{thm1.2} completely.

Note that Theorem \ref{thm1.2} says that the manifold has at most linear volume growth.
We remark that Calabi \cite{Cala75} and Yau \cite{Y76} independently proved that every noncompact manifold with nonnegative Ricci curvature has at least linear volume growth.

In \cite{Y92}, Yau proposed the following question.
\begin{prob}\label{prob1.3}
Let $(M^n,g)$ be an $n$-dimensional complete non-compact Riemannian manifold with nonnegative Ricci curvature. Then for any fixed point $p\in M$, do we have
\begin{equation}
\limsup_{r\rightarrow\infty}r^{2-n}\int_{B(p,r)}\mathrm{Sc}_g<\infty?
\end{equation}
\end{prob}
It is easy to see that a positive answer of Question \ref{prob1.3} will give an answer to Conjecture \ref{prob1.1}.

In \cite{Z22}, Zhu proved the correctness of Question \ref{prob1.3} for 3-dimensional manifolds with a pole (recall that $M$ has a pole $p\in M$ means that the exponential map at $p$ is a diffeomorphism) and $\mathrm{Ric}_g\geq 0$.
In \cite{Xu2306}, Xu proved the following sharp result:
\begin{thm}\label{thm1.4}
Let $(M,g)$ be a 3-dimensional Riemannian manifold with a pole. Suppose $\mathrm{Ric}_g\geq 0$, then we have
\begin{align}
\lim_{r\rightarrow \infty}\frac{\int_{B(p,r)}\mathrm{Sc}_{g}}{r}=8\pi(1-\lim_{r\rightarrow\infty}\frac{\mathrm{Vol}(B(p,r))}{\omega_3 r^{3}}),
\end{align}
where $\omega_3$ is the volume of unit ball in $\mathbb{R}^3$.
\end{thm}

The aim of the present paper consists of two parts.
Firstly, in Theorem \ref{thm1.2}, the universal constant $C$ is not explicit.
Motivated by Theorem \ref{thm1.4}, we want to find an optimal volume ratio near infinity.
Secondly, we want to explore the influence of positive scalar curvature on the geometry on a larger class of manifolds which may have negative Ricci curvature somewhere.

The following theorem is our main result.

\begin{thm}\label{thm-opti-vol}
Let $(M^3,g)$ be a $3$-dimensional complete noncompact Riemannian manifold with $k$ ends and finite first Betti number.
Suppose $S_g\geq2$ outside a compact set, and for some fixed pointed $p\in M$, there is a positive continuous function $f:[0,\infty)\rightarrow (0,\infty)$ such that $\displaystyle\lim_{r\rightarrow\infty}f(r)=0$ and
$$\mathrm{Ric}_g(x)\geq-f(r(x)),$$
where $r(x)=d(x,p)$.
Then the following conclusions were justified.

(1)
There exists $C>0$ such that
\begin{equation}\label{mthm2}
\limsup_{r\rightarrow\infty}\frac{\mathrm{Vol}(B(p,r))}{r}\leq C.
\end{equation}

(2)If in addition $\displaystyle f(r)$ is a non-increasing function and
\begin{align}\label{decay}
\int_{0}^{\infty}rf(r)dr<\infty,
\end{align}
then
\begin{equation}\label{opti-vol}
\limsup_{r\rightarrow\infty}\frac{\mathrm{Vol}(B(p, r))}{r}\leq4k\pi.
\end{equation}
\end{thm}



\begin{rem}
Recall that any manifolds with nonnegative Ricci curvature has at most 2 ends.
In addition, $3$-dimensional manifolds with nonnegative Ricci curvature have been classified by Liu in \cite{L13}, and it is known that every such manifold must have finite first Betti number.
Therefore, Theorem \ref{thm-opti-vol} holds automatically on 3-manifolds with nonnegative Ricci curvature and positive scalar curvature.
In this case, (2) in Theorem \ref{thm-opti-vol} gives a sharp upper bound.
If $k=2$, then according to Cheeger-Gromoll's splitting theorem, $M$ is isometric to $\mathbb{R}\times S$ with $S$ a compact manifold with sectional curvature $\geq 1$, and it is easy to see that the equality in (\ref{opti-vol}) holds only when $S$ is isometric to $\mathbb{S}^2(1)$.
If $k=1$, then the equality in (\ref{thm-opti-vol}) holds on $M=(\mathbb{R}\times\mathbb{S}^2(1))/\mathbb{Z}_2$, where $\mathbb{Z}_2=\langle1,\sigma\rangle$, and $\sigma:\mathbb{R}\times\mathbb{S}^2(1)\rightarrow\mathbb{R}\times\mathbb{S}^2(1)$, $\sigma(r,s)=(-r,\tau(s))$, where $\tau$ is the antipodal map on $\mathbb{S}^2(1)$.
\end{rem}

(1) of Theorem \ref{thm-opti-vol} can be proved by the same ideas in \cite{CLS23}, while (2) of Theorem \ref{thm-opti-vol} needs some new observations.
In the proof of (2) in Theorem \ref{thm-opti-vol}, we first prove the following proposition.

\begin{prop}\label{b}
Let $(M^n,g)$ be an $n$-dimensional complete noncompact Riemannian manifold.
Given $p\in M$, denote by $r(x)=d(x,p)$.
Suppose there is a non-increasing continuous function $f:[0,\infty)\rightarrow (0,\infty)$ satisfying (\ref{decay}), and
\begin{equation} \label{Ric-decay}
\mathrm{Ric}_g(x)\geq-f(r(x))
\end{equation}
holds on $M$.
Let $\gamma:[0,\infty)\rightarrow M$ be a ray and $b$ be the Busemann function associated to $\gamma$.
Then there exists a constant $C_1>0$ such that
\begin{equation}\label{1.3}
\mathrm{Vol}(b^{-1}([k,k+1]))>C_1
\end{equation}
for any sufficiently large $k$.
In particular, $M$ has at least linear volume growth, i.e. there exists a constant $c>0$ such that
\begin{equation}\label{1.4}
\liminf_{r\rightarrow\infty}\frac{\mathrm{Vol}(B_r(p))}{r}>c.
\end{equation}
\end{prop}

Manifolds with the Ricci curvature lower bound (\ref{Ric-decay}), where $f$ satisfies (\ref{decay}), have been studied in many previous works such as \cite{AG90}, \cite{LiT92}, \cite{Z94}, \cite{MW22} etc.
In particular, according to Lemma 3.1 in \cite{LiT92}, such a manifold has at least linear volume growth, which generalizes the famous result of Calabi \cite{Cala75} and Yau \cite{Y76}.
In Proposition \ref{b}, we utilize the method presented by Sormani in \cite{S1} to prove (\ref{1.3}), and this gives another proof of the Calabi-Yau type property.

\begin{rem}
By the Calabi-Yau type property in \cite{Cala75} and \cite{Y76}, it is interesting to study manifolds (or metric measure spaces) with nonnegative Ricci curvature and linear volume growth.
There are many papers concerning this class of manifolds (or metric measure spaces), see \cite{AP23}, \cite{H18}, \cite{H20}, \cite{S1}, \cite{Sor99}, \cite{Zxy23} etc.
\end{rem}

Now we sketch a proof of (\ref{opti-vol}) under the assumptions that $M$ has nonnegative Ricci curvature and $\mathrm{Sc}_g\geq2$ \textbf{everywhere}, and $k=1$.
Firstly, we choose a sequence of points $p_{i}$ lying on a ray with $p_{i}$ going to infinity.
Then by Cheeger-Colding's theory, up to a subsequence, $(M,d_g,p_i)$ converges to a Ricci limit space $(X,d_X)$, which is isometric to $(\mathbb{R}\times Y, d_{\mathrm{Eucl}}\otimes d_{Y})$
for some metric space $(Y,d_Y)$.
We need to understand more properties of $Y$.
The properties of noncollapsed Ricci limit spaces coming from manifolds with a uniformly positive scalar curvature lower bound were studied in the recent papers \cite{WXZZ22} \cite{ZZ23} etc.
The limit spaces which are studied by \cite{WXZZ22} \cite{ZZ23} are assumed to be noncollapsed everywhere, i.e. the volumes of \textbf{any} unit ball have a uniformly positive lower bound.
In addition, according to the proof in \cite{ZZ23}, which uses the theory of Ricci flow coming from 3-manifolds with possibly unbounded curvature (see \cite{ST17} and \cite{BCW17} etc.), it can be shown that if the above $X$ is noncollapsed everywhere, then the above $(Y,d_Y,\mathcal{H}^{2})$ is in fact a non-collapsed $\mathrm{RCD}(1,2)$ space (in the sense of \cite{GigPhi15}), then we have an optimal upper bound of $\mathcal{H}^{2}(Y)$, which is a key in the proof of (\ref{opti-vol}).

In the proof of (\ref{opti-vol}), we prove the following proposition:
\begin{prop}\label{lem-non-col}
Under the assumption of (2) in Theorem \ref{thm-opti-vol}, there is $\delta>0$ such that for all $q\in M$,
\begin{equation}\label{1.9}
\mathrm{Vol}(B(q,1))>\delta.
\end{equation}
\end{prop}

\begin{rem}
Proposition \ref{lem-non-col} will be proved in Section \ref{sec-5}. From the proof of Proposition \ref{lem-non-col}, we can see that, if a higher dimensional manifold $M^{n}$ has a Ricci curvature lower bound as in Proposition \ref{b}, and assume in addition that for an end $E$ of $M$, there exists $C>0$ such that
$\mathrm{diam}(E\cap \partial B(p,R))\leq C$
for every sufficiently large $R$,
then there exists $\delta>0$ such that
$\mathrm{Vol}(B(q,1))>\delta$
holds for every $q\in E$.

\end{rem}

Proposition \ref{lem-non-col} enables us to apply the theory of Ricci flow, as what is done in \cite{ZZ23}.
Our proof of (\ref{opti-vol}) under the general assumptions in Theorem \ref{thm-opti-vol} still uses Ricci flow.
Because the lower bound of Ricci curvature may be negative some where and the positive lower bound of scalar curvature is not global, we must make some adjustments to the proof in \cite{ZZ23}.

\begin{rem}
The results of this paper are part of the second author's thesis.
During the process of writing this paper, Wei, Xu and Zhang have posted their paper \cite{WXZ24} on arXiv.
In the independent work \cite{WXZ24}, the sharp volume ratio near infinity has also been obtained on manifolds with $\mathrm{Ric}\geq0$ and $\mathrm{Sc}\geq2$.
The methods used in \cite{WXZ24} and our paper are different.
In addition, the corresponding rigidity problem and examples are discussed in \cite{WXZ24}, while we consider manifolds with asymptotically nonnegative Ricci curvature in our paper.

After we posted the first version of our preprint onto arXiv, we find a paper (\cite{AX24}) of Antonelli and Xu, which is posted onto arXiv on the same day with ours.
In this independent work, Antonelli and Xu have obtained many interesting results, one of which also considers the sharp asymptotic volume ratio near infinity on manifolds with dimension at most $5$, see Remark 7 of \cite{AX24}.
We remark that in \cite{AX24}, to obtain the sharp asymptotic volume ratio, the assumption on nonnegative Ricci curvature is required in every dimension at most $5$, and in dimension $4$ and $5$, the assumption on the uniformly positive lower bound of scalar curvature is replaced by a uniformly positive lower bound of the biRicci curvature (or some kind of its generalization to the spectral sense).





\end{rem}
\vspace*{12pt}

In the following we give the outline of this paper.
Section \ref{sec-2} contains some preliminaries on distance functions and Busemann functions.
In Section \ref{sec-3}, we prove Proposition \ref{b}.
In Section \ref{sec-4}, we give the proof of part (1) of Theorem \ref{thm-opti-vol}.
In Section \ref{sec-5}, we prove Proposition \ref{lem-non-col}.
In Section \ref{sec-6}, we finish the proof of part (2) of Theorem \ref{thm-opti-vol}.

\vspace*{12pt}

\noindent\textbf{Acknowledgments.}

The authors would like to thank Profs. Hongzhi Huang and Shaochuang Huang for discussions.
The authors are partially supported by National Key R\&D Program of China(2021YFA1002100) and National
Natural Science Foundation of China (Nos. 12271531 and 12025109).

\section{Preliminaries}\label{sec-2}
The manifold we are considering in this section is always complete and noncompact.
All geodesics are parametrized by arclength.

A geodesic $\gamma:(-\infty,\infty)\rightarrow M$ is called a ray  if for any interval $[a,b]\subset(-\infty,\infty)$, $\gamma|_{[a,b]}$ is a minimizing geodesic.
A ray is a half-line $\gamma:[0,\infty)\rightarrow M$.

Given a ray $\gamma$, for any $t\geq0$, denote by $b_{t}(x):=t-d(x,\gamma(t))$.
It is easy to check that, given any $x\in M$, the function $t\mapsto b_{t}(x)$ is non-decreasing, and $b_{t}(x)\leq d(x,p)$ for all $t\geq0$, hence
\begin{equation}
b_{\gamma}(x)=\lim_{t\rightarrow\infty}b_{t}(x)
\end{equation}
is well-defined, which is called the Busemann function associated to $\gamma$, and is abbreviated as $b(x)$ if there is no confusion on $\gamma$.

A ray $\gamma_x$ emanating from $x$ is called a Busemann ray associated with $\gamma$ if it is the limit of a sequence of minimal geodesics, $\sigma_i$, from $x$ to $\gamma(R_i)$ in the following sense,
\begin{equation}
\gamma'_x(0)=\lim_{R_i\rightarrow\infty}\sigma'_i(0).
\end{equation}
We reparametrize $\gamma_x$ by arc length so that $\gamma_x(b(x))=x$.



Let $K$ be a compact set contained in $M$. Then we define
\begin{equation}
\Omega(K)=\{x|\exists z\in K,\exists t\geq b(z),\exists \gamma_z\; \text{such that } x=\gamma_z(t) \}.
\end{equation}

Furthermore, let $\Omega_R(K)=\Omega(K)\cap b^{-1}(R)$.
It is not hard to check that $\Omega(K)$ is closed and $\Omega_R(K)$ is compact.
The following lemma shows that the Busemann function can be regarded as a distance function.

\begin{lem}[Lemma 6 in \cite{S1}]\label{sormani6}
Fix $r<R$ and $K$ be a subset of $b^{-1}((-\infty,r))$. Then
\begin{equation}
d(x,b^{-1}(R))=R-b(x),\; \forall x\in b^{-1}((-\infty,R]).
\end{equation}
Furthermore, if $z\in \Omega(K)\cap b^{-1}([r,R])$, then there is only one point, $y\in b^{-1}(R)$, such that $d(z,y)=d(z,b^{-1}(R))$.
Thus $y=\gamma_z(R)$ and
\begin{equation}
d(z,\Omega_R(K))=R-b(z),\; \forall z\in \Omega(K)\cap b^{-1}([r,R]).
\end{equation}
\end{lem}

Now we are proving that under certain conditions, the distance function and the Busemann function can control each other.

\begin{lem}\label{br}
Let $\gamma$ be a ray emitting from $p\in M$ $(\gamma(0)=p)$, and let $b(x)$ be the Busemann function associated to $\gamma$.
Then
\\(a)$b^{-1}(R)\subset \overline{M\setminus B(p,R)},\;\forall R>0$;
\\(b)Let $E$ be an end of $M$ such that $\gamma(t)\in E$ for every sufficiently large $t$.
Then for every sufficiently large $R>0$, $b^{-1}(R)\subset E$.
In addition, if there exists a positive constant $C>0$ such that $\mathrm{diam}(E\cap\partial B(p,r))\leq C$ for every sufficiently large $r>0$, then for any sufficiently large $R$ it holds that
\begin{equation}
b^{-1}(R)\subset\overline{E\cap B(p,R+C)\setminus B(p,R)},
\end{equation}
and in particular,
\begin{equation}
\mathrm{diam}(b^{-1}(R)) \leq 3C.
\end{equation}
\end{lem}
\begin{proof}
For (a), suppose there exists a point $x\in b^{-1}(R)\cap B(p,R)$, then $d(p,x)<R$. Since $b(p)=0$, by Lemma \ref{sormani6} we know that
\begin{equation}
d(p, b^{-1}(R))=R,
\end{equation}
which leads to
\begin{equation}
R>d(p,x) \geq d(p,b^{-1}(R))=R,
\end{equation}
resulting in a contradiction.

In the following we prove (b).

\textbf{Claim 1.} Given any $R$. Then any $x,y\in b^{-1}(R)$ can be connected by a curve contained in $b^{-1}([R,\infty))$.

\textbf{Proof of Claim 1.} We only need to prove that any $x\in b^{-1}(R)$ can be connected to $\gamma(R+1)$ by a curve $\sigma\subset b^{-1}([R,\infty))$.
Firstly, $\gamma_{x}(t)\subset b^{-1}([R,\infty))$ connects $x$ and $y=\gamma_{x}(R+1)$.
Since $R+1=b(y)=\lim_{t\rightarrow\infty}b_{t}(y)$, we choose $\tilde{t}$ sufficiently large such that $b_{\tilde{t}}(y)=\tilde{t}-d(y,\gamma(\tilde{t}))>R+\frac{1}{2}$.
Let $\tilde{\sigma}:[0,d(y,\gamma(\tilde{t}))]\rightarrow M$ be a geodesic connecting $y$ to $\gamma(\tilde{t})$.
Then for any $r\in[0,d(y,\gamma(\tilde{t}))]$, we have $b(\tilde{\sigma}(r))\geq b_{\tilde{t}}(\tilde{\sigma}(r))=b_{\tilde{t}}(y)+r>R+\frac{1}{2}$.
Hence $\tilde{\sigma}\subset b^{-1}([R,\infty))$.
Finally, $\gamma$ itself is a curve connecting $\gamma(\tilde{t})$ and $\gamma(R+1)$.
Glue the above curves to form a desired curve and this completes the proof of Claim 1.

Now we assume $R$ is large enough so that a unbounded connected component of $M\setminus B(p,R)$, denoted by $N$, is contained in $E$.
According to (a) and Claim 1, every point in $b^{-1}(R)$ can be connected to $\gamma(R)$ by a curve contained in $M\setminus B(p,R)$.
Since $\gamma(R)\in N$, we have $b^{-1}(R)\subset N\subset E$.

\textbf{Claim 2. } Suppose $R>0$ is sufficiently large, then $E\cap \partial B(p,R)\subset b^{-1}([R-C,R])$.

\textbf{Proof of Claim 2.} By (a), we have $\partial B(p,R)\subset b^{-1}((-\infty,R])$.
Suppose the claim does not hold, then there exists $S<R-C$ such that $E\cap \partial B(p,R) \cap b^{-1}(S)\neq\varnothing$.
Take $y\in E\cap\partial B(p,R)\cap b^{-1}(S)$.
Since $\mathrm{diam}(E\cap\partial B(p,R))\leq C$ and $\gamma(R) \in \partial B(p,R) \cap b^{-1}(R)$, we know that
\begin{equation}
d(y,b^{-1}(R))\leq d(y,\gamma(R))\leq C.
\end{equation}
On the other hand, by Lemma \ref{sormani6} we know that $d(y, b^{-1}(R))=R-S>C$, which leads to a contradiction. This completes the proof of Claim 2.

Now, let's prove that $\displaystyle\sup_{x\in b^{-1}(R)}d(p,x)\leq R+C$.
Suppose the conclusion does not hold, then there exists a point $z\in b^{-1}(R)\subset E$ and $\varepsilon>0$ such that $d(p,z)=R+C+\varepsilon$.
Hence $E\cap\partial B(p,R+C+\varepsilon)\cap b^{-1}(R)\neq\varnothing$.
On the other hand, according to Claim 2, we have $E\cap\partial B(p,R+C+\varepsilon)\subset b^{-1}([R+\varepsilon,R+C+\varepsilon])$, which is a contradiction.
Therefore, $b^{-1}(R) \subset E\cap \overline{B(p,R+C)\setminus B(p,R)}$. Hence,
\begin{align}
  \mathrm{diam}(b^{-1}(R))&\leq \mathrm{diam}(E\cap\overline{B(p,R+C) \setminus B(p,R)})\\
  & \leq C+\mathrm{diam}(E\cap\partial B(p,R)) + \mathrm{diam}(E\cap\partial B(p,R+C)) \nonumber\\
  & \leq 3C. \nonumber
\end{align}
The proof is completed.
\end{proof}

Given any $\delta>0$ and any $r_1<r_2<R$, let
\begin{equation}
S_{\delta,r_1,r_2}=S_{\delta,r_1,r_2}(\Omega_R(K))
\end{equation}
be the set of points $x$ with $d(x,\Omega_R(K))\in[R-r_2,R-r_1]$ such that there exists a minimal geodesic $\sigma$ from $\Omega_R(K)$ to $x$ with
\begin{align}
\mathrm{Length}(\sigma)=d(x,\Omega_R(K))\quad\text{ and }\quad g(\sigma'(0),-\nabla b)\geq 1-\delta.
\end{align}

In the following we collect some notions and results from \cite{S1}.

\begin{lem}[Lemma 8 and Corollary 9 in \cite{S1}]
Fix $R>r_2>r_1$, $\delta>0$ and $K$ be a compact subset of $b^{-1}((-\infty,r_1))$. Then
\begin{equation}
\bigcap_{\delta>0}S_{\delta,r_1,r_2}(\Omega_R(K))=\Omega(K)\cap b^{-1}([r_1,r_2])
\end{equation}
and
\begin{equation}\label{Sb}
\lim_{\delta\rightarrow0}\mathrm{Vol}(S_{\delta,r_1,r_2}(\Omega_R(K)))=\mathrm{Vol}(\Omega(K)\cap b^{-1}([r_1,r_2])).
\end{equation}
\end{lem}

\begin{lem}[Lemma 10 in \cite{S1}]\label{lem2.3}
Let $K\subset b^{-1}((-\infty,r_1])$ and $R>r_3>r_2>r_1$. Given any $\delta>0$, there exists
\begin{equation}
h_1(\delta)=h_1(\delta,r_1,r_2,r_3,R,K,M)<\delta
\end{equation}
such that for any $x\in S_{h_{1}(\delta),r_2,r_3}(\Omega_R(K))$, every minimal geodesic, $\sigma$, from $\sigma(0)\in\Omega_R(K)$ to $x$ with $d(x,\sigma(0))=d(x,\Omega_R(K))$, it holds
\begin{equation}
g(\sigma'(0),-\nabla b)\geq 1-\delta.
\end{equation}
\end{lem}


Given any $\varepsilon>0$, there exists a set
\begin{equation}
\Omega_R^{\varepsilon}(K)=\{p_1,\cdots, p_{N_{\varepsilon}}\}\subset\Omega_R(K)
\end{equation}
such that the tubular neighborhood, $T_{\varepsilon}(\Omega_R^{\varepsilon}(K))$, contains $\Omega_R(K)$.
In this paper, for a set $X\subset M$ and $0\leq a<b$, we use the notations
$$T_{a,b}(X):=\{y|a\leq d(y,X)\leq b\},$$
and
$$T_{b}(X):=T_{0,b}(X).$$

Given $\Omega_R^{\varepsilon}(K)$ and $p\in \Omega_R^{\varepsilon}(K)$, the star-shaped set $V_{p,\varepsilon}$ and the star-shaped wedge set $U_{p,\delta}$ are defined to be
\begin{equation}
V_{p,\varepsilon}=\{x|d(x,p)<d(x,q),\;\forall q\in \Omega_R^{\varepsilon}(K)\; \text{such that } q\neq p\},
\end{equation}
\begin{equation}
U_{p,\delta}=\{x|\exists \text{ a minimal geodesic } \sigma \text{ from }p \text{ to } x \text{ such that } g(\sigma'(0),-\nabla b)\geq 1-\delta\}.
\end{equation}
In addition, denote by
\begin{equation}\label{Uvd}
U_{\varepsilon,\delta}=\bigcup_{i=1}^{N_{\varepsilon}}(U_{p_i,\delta}\cap V_{p_i,\varepsilon}).
\end{equation}

\begin{lem}[Lemma 14 in \cite{S1}]\label{lem2.4}
Fix $r_1<R$ and the compact set, $K\subset b^{-1}((-\infty,r_1])$.
For all $\varepsilon>0$, there exists $h_2(\varepsilon,r_1,R,K,M)>0$ such that, for all $p\in\Omega^{\varepsilon}_R(K)$ and all $\delta<h_2$ we have
\begin{equation}
d(\gamma_x(R),p)<\varepsilon,\quad\forall x\in U_{p,\delta}\cap B(p,R-r_1-\varepsilon).
\end{equation}
Thus
\begin{equation}
b(x)\in[R-d(x,p),R-d(x,p)+\varepsilon],\quad\forall x\in U_{p,\delta}\cap B(p,R-r_1-\varepsilon).
\end{equation}
\end{lem}

\begin{lem}[Lemma 15 in \cite{S1}]\label{bT}
Fix $r_1<r_2<R$ and $\delta>0$. Let $\varepsilon_0>0$, $K\subset b^{-1}((-\infty,r_1])$, there exists $h_3(M,\delta,r_1,r_2,R)>0$ sufficiently small such that for all $\varepsilon<\min\{\varepsilon_0,h_3\}$,
\begin{equation}
\Omega(K)\cap b^{-1}([r_1,r_2])\subset T_{R-r_2-\varepsilon_0,R-r_1+\varepsilon_0}(\Omega_R^{\varepsilon}(K))\cap\bar{U}_{\varepsilon,\delta}.
\end{equation}
where $\bar{U}$ is the closure of $U$.
So
\begin{equation}
\mathrm{Vol}(\Omega(K)\cap b^{-1}([r_1,r_2]))\leq \mathrm{Vol}(T_{R-r_2-\varepsilon_0,R-r_1+\varepsilon_0}(\Omega_R^{\varepsilon}(K))\cap U_{\varepsilon,\delta}).
\end{equation}
\end{lem}

\begin{lem}[Lemma 16 in \cite{S1}]\label{TS}
Fix $r_1<r_2<r_3<R$ and $K\subset b^{-1}((-\infty,r_1])$. Given any $\delta>0$, let $h_1(\delta)<\delta$ be the constant defined in Lemma \ref{lem2.3}.
Then there exists $h_4=h_4(M,h_1(\delta),R,r_1,r_2,r_3)>0$ such that given any $\varepsilon_0>0$, for all $\varepsilon<\min\{h_4,\varepsilon_0\}$, we have
\begin{equation}
T_{R-r_3+\varepsilon_{0},R-r_2-\varepsilon_0}(\Omega_R^{\varepsilon}(K))\cap U_{\varepsilon,h_1(\delta)}\subset S_{2\delta,r_2,r_3}(\Omega_R(K)).
\end{equation}
\end{lem}

\section{Proof of Proposition \ref{b}}\label{sec-3}
\begin{prop}\label{3.1}
Let $(M^n,g)$ be an $n$-dimensional complete noncompact Riemannian manifold with a given ray $\gamma$ emitting from $q$ and its associated  Busemann function, $b=b_{\gamma}$ with $b(\gamma(t))=t$.
Suppose there is a non-increasing continuous function $f:[0,\infty)\rightarrow (0,\infty)$ satisfying (\ref{decay}), and
\begin{equation} \label{Ric-decay-Bf}
\mathrm{Ric}_g(x)\geq-f(b(x)),\quad \forall x\in b^{-1}([0,\infty)).
\end{equation}
Let $K\subset b^{-1}((0,\infty))$ be a compact set, then there exists a constant $C_1>0$ such that
\begin{equation}\label{3.3}
\mathrm{Vol}(\Omega(K)\cap b^{-1}([k,k+1]))>C_1
\end{equation}
for any sufficiently large $k$.
\end{prop}
\begin{proof}
Firstly, we define a positive function $k(r) =\int_{r}^{\infty}f(s)ds$.
Since $f>0$ satisfies (\ref{decay}), it is easy to see that $k(r)$ is well defined and $\int_{0}^{\infty}k(s)ds=\int_{0}^{\infty}sf(s)ds<\infty$.

Let $u_{R}$ be the solution of
\begin{align}
\varphi'(r)=\varphi^{2}(r)-f(r)
\end{align}
on $[0, R]$ with boundary condition $u_{R}(R) = 0$.
Since $\varphi= 0$ is a subsolution and $k$ is a supersolution, it can be shown that $0\leq u_{R} \leq k$ on $[0, R]$.
Hence $u := \lim_{R\rightarrow \infty}u_{R}$ exists and satisfies
\begin{align}
u'(r)=u^{2}(r)-f(r)
\end{align}
on $[0,\infty)$, and $0\leq u(r) \leq k(r)$ for every $r > 0$.
In particular,
\begin{align}
\int_{0}^{\infty}u(r)dr\leq \int_{0}^{\infty}k(r)dr<\infty.
\end{align}

Now we define
\begin{align}\label{z-function}
z(r) = \exp(- \int_{0}^{r}u(s)ds),
\end{align}
then it is easy to see that $z(0) = 1$, $z(r) > 0$, $z(r)$ is non-increasing and
\begin{align}
z''(r)= f(r)z(r).
\end{align}

\begin{lem}\label{lem3.1}
Let $f: [0,\infty)\rightarrow(0,\infty)$ be a continuous non-increasing function such that (\ref{decay}) holds.
For any $R>0$, let $g$ be the solution of the equation
\begin{align}
g''(r) = f(R-r)g(r)
\end{align}
on $[0, R]$ with initial conditions $g(0) = 0$ and $g'(0) = 1$.
Then
\begin{align}\label{3.9}
z(R-r)-\frac{R-r}{R}z(R) \leq \frac{g(r)}{g(R)}\leq z(R-r)
\end{align}
where $z$ is given in (\ref{z-function}).
\end{lem}

\begin{proof}[Proof of Lemma \ref{lem3.1}.]
Define $h(r) = \frac{g(R-r)}{g(R)}$.
Then $h\geq0$ satisfies
\begin{align}
h''(r) = f(r)h(r)
\end{align}
on $[0,R]$ with boundary conditions $h(0)=1$ and $h(R) = 0$.
By the maximum principle, we have
\begin{align}\label{3.11}
z(r)\geq h(r)
\end{align}
on $[0, R]$.
Hence $z-h$ is convex.
In particular, for any $r\in[0, R]$, we have
\begin{align}
z(r)-h(r)\leq \frac{r}{R}(z(R)-h(R)) + (1-\frac{r}{R})(z(0)-h(0))\leq \frac{r}{R}z(R),
\end{align}
and hence
\begin{align}\label{3.13}
h(r)\geq z(r)-\frac{r}{R}z(R).
\end{align}

(\ref{3.9}) follows from (\ref{3.11}) and (\ref{3.13}).
\end{proof}

Proposition \ref{3.1} can be proved provided the following claim is justified.

\textbf{Claim 1:} There exists some $k_{0}>0$ such that for any $k>k_{0}$, it holds
\begin{equation}\label{3.2}
\mathrm{Vol}(\Omega(K)\cap b^{-1}([k+3,k+4]))\geq \bigl(\frac{z(k+4)}{z(k+1)}\bigr)^{n-1}\mathrm{Vol}(\Omega(K)\cap b^{-1}([k+1,k+2])).
\end{equation}

\begin{proof}[Proof of Claim 1:]
Let
\begin{equation}
r_4>r_3>r_2>r_1> \sup_{x\in K}b^{-1}(x).
\end{equation}
Take $R>r_{4}$ sufficiently large such that $z(r_{4})-\frac{r_{4}}{R}>0$, where $z$ is given in (\ref{z-function}).
Take $\displaystyle\varepsilon_0<\frac{1}{10}\min\{r_2-r_1,r_4-r_3,R-r_4\}$.
By Lemma \ref{lem2.4}, there exists $h_2:=h_{2}(\varepsilon_{0},r_{1},R,K,M)$, such that if $p\in \Omega_R^{\varepsilon}(K)$ and $0<\delta<h_2$, then for all $0<h<\delta$ and $x\in U_{p,h}\cap B(p,R-r_1-\varepsilon_0)$, we have
\begin{equation}\label{bR}
 b(x)\in[R-d(x,p),R-d(x,p)+\varepsilon_0].
\end{equation}

Take $h_1:=h_{1}(\delta,r_{1},r_{3}-\varepsilon_{0},r_{4}+\varepsilon_{0},R,K,M)<\delta$ from Lemma \ref{lem2.3}.
Then (\ref{bR}) holds for every $x\in U_{p,h_1}\cap B(p,R-r_1-\varepsilon_0)$, and hence
\begin{equation}\label{3.17}
\mathrm{Ric}_g(x)\geq -f(R-d(x,p))
\end{equation}
holds for these $x$.

For every fixed large $R>0$, let $J_{R}:[0,R]\rightarrow\mathbb{R}^{+}$, be the solution of
\begin{align}
J_{R}''(t)=f(R-t)J_{R}(t)
\end{align}
with the initial conditions $J_{R}(0) = 0$, $J_{R}'(0) = 1$.

Let $\mathrm{M}^{J_{R}}=[0,R]\times_{J_{R}} \mathbb{S}^{n-1}$ be the warped product spaces equipped with the metric $g_R=dt^2+J_{R}^2g_{\mathbb{S}^{n-1}}$, where $g_{\mathbb{S}^{n-1}}$ is the standard metric on $\mathbb{S}^{n-1}$.
Then at any $z=(t,s)\in \mathrm{M}^{J_{R}}$, we have
\begin{align}
\mathrm{Ric}_{g_0}(\partial t,\partial t)=-(n-1)\frac{J_{R}''(t)}{J_{R}(t)} = -(n-1)f(R-t),
\end{align}
\begin{align}
\mathrm{Ric}_{g_0}(v,v)=(n-2)\frac{1-J_{R}'^2}{J^2_{R}}-\frac{J_{R}''}{J_{R}}\leq -f(R-t),
\end{align}
where $v$ is any unit vector orthogonal to $\partial t$.
Since (\ref{3.17}) holds, for every $z\in \mathrm{M}^{J_{R}}$ and $y\in U_{p,h_1}\cap B(p,R-r_1-\varepsilon_{0})$ with $d(y,p)=d_{g_{R}}(z,0)$, we have
\begin{equation}\label{3.12}
\mathrm{Ric}_g(y)\geq \mathrm{Ric}_{g_R}(z).
\end{equation}

Choose $h_3=h_3(M,h_{1},r_{1}+2\varepsilon_{0},r_{2}-\varepsilon_{0},R)$ and $h_4=h_4(M,h_{1},R, r_{1}, r_{3}-\varepsilon_{0},r_{4}+\varepsilon_{0})$ from Lemmas \ref{bT} and \ref{TS} respectively.
Let
\begin{equation}\label{07}
\varepsilon<\min\{h_3,h_4,\varepsilon_0\}.
\end{equation}
Then we have
\begin{equation}\label{3.222}
\frac{\mathrm{Vol}(S_{2\delta,r_3-\varepsilon_0,r_4+\varepsilon_0}(\Omega_R(K)))}{\mathrm{Vol}(\Omega(K)\cap b^{-1}([r_1+2\varepsilon_0,r_2-\varepsilon_0]))}
\geq \frac{\mathrm{Vol}(T_{R-r_4,R-r_3}(\Omega_R^{\varepsilon}(K)))\cap U_{\varepsilon,h_1})}{\mathrm{Vol}(T_{R-r_2,R-r_1-\varepsilon_0}(\Omega_R^{\varepsilon}(K)))\cap U_{\varepsilon,h_1})}.
\end{equation}

For $p\in \Omega_R^{\varepsilon}(K)$, denote by
\begin{equation}\label{Vp}
V_p=V_{p,\varepsilon}\cap U_{p,h_1},
\end{equation}
For any $r<s$, denote by $\mathrm{Ann}_p(r,s)=B(p,s)\setminus B(p,r)$.
Then for $s<R-r_1-\varepsilon_0$, we have the disjoint union
\begin{equation}
T_{r,s}(\Omega_R^{\varepsilon}(K))\cap U_{\varepsilon,h_1}=\bigcup_{p\in \Omega_R^{\varepsilon}(K)}(\mathrm{Ann}_p(r,s)\cap V_p).
\end{equation}
Hence
\begin{equation}
\mathrm{Vol}(T_{r,s}(\Omega_R^{\varepsilon}(K))\cap U_{\varepsilon,h_1})=\sum_{p\in \Omega_R^{\varepsilon}(K)}\mathrm{Vol}(\mathrm{Ann}_p(r,s)\cap V_p).
\end{equation}

Since (\ref{3.12}) holds, applying the volume comparison theorem to starlike sets as in \cite{C90}, we obtain
\begin{equation}
\frac{\mathrm{Vol}(\mathrm{Ann}_p(R-r_4,R-r_3)\cap V_p)}{\mathrm{Vol}(\mathrm{Ann}_p(R-r_2,R-r_1-\varepsilon_0)\cap V_p)}\geq \frac{V(r_3,r_4,R)}{V(r_1+\varepsilon_0,r_2,R)},
\end{equation}
where \begin{align}\label{V=}
V(s_1,s_2,R)=\int_{R-s_2}^{R-s_1}J_{R}^{n-1}(t)dt.
\end{align}

Thus,
\begin{align}\label{TV}
&\frac{\mathrm{Vol}(T_{R-r_4,R-r_3}(\Omega_R^{\varepsilon}(K))\cap U_{\varepsilon,h_1})}{\mathrm{Vol}(T_{R-r_2,R-r_1-\varepsilon_0}(\Omega_R^{\varepsilon}(K))\cap U_{\varepsilon,h_1})} \\
=&\frac{\sum_{p\in \Omega_R^{\varepsilon}(K)}\mathrm{Vol}(\mathrm{Ann}_p(R-r_4,R-r_3)\cap V_p)}{\sum_{p\in \Omega_R^{\varepsilon}(K)}\mathrm{Vol}(\mathrm{Ann}_p(R-r_2,R-r_1-\varepsilon_0)\cap V_p)}\nonumber\\
   \geq& \frac{\sum_{p\in \Omega_R^{\varepsilon}(K)}\mathrm{Vol}(\mathrm{Ann}_p(R-r_2,R-r_1-\varepsilon_0)\cap V_p)\frac{V(r_3,r_4,R)}{V(r_1+\varepsilon_0,r_2,R)}}{\sum_{p\in \Omega_R^{\varepsilon}(K)}\mathrm{Vol}(\mathrm{Ann}_p(R-r_2,R-r_1-\varepsilon_0)\cap V_p)} \nonumber\\
   =&\frac{V(r_3,r_4,R)}{V(r_1+\varepsilon_0,r_2,R)}.\nonumber
\end{align}

Combining (\ref{3.222}) with (\ref{TV}), we have
\begin{equation}\label{sbv}
\frac{\mathrm{Vol}(S_{2\delta,r_3-\varepsilon_0,r_4+\varepsilon_0}(\Omega_R(K)))}{\mathrm{Vol}(\Omega(K)\cap b^{-1}([r_1+2\varepsilon_0,r_2-\varepsilon_0]))}
 \geq \frac{V(r_3,r_4,R)}{V(r_1+\varepsilon_0,r_2,R)}.
\end{equation}
Note that (\ref{sbv}) does not depend on $\delta$ or $\varepsilon$.
Taking the limit as $\delta\rightarrow0$ and applying (\ref{Sb}), we obtain
\begin{equation}\label{bV}
 \frac{\mathrm{Vol}(\Omega(K)\cap b^{-1}([r_{3}-\varepsilon_0,r_{4}+\varepsilon_0]))}{\mathrm{Vol}(\Omega(K)\cap b^{-1}([r_{1}+2\varepsilon_0,r_{2}-\varepsilon_0]))}\geq \frac{V(r_3,r_4,R)}{V(r_1+\varepsilon_0,r_2,R)}.
\end{equation}

On the other hand,
\begin{align}
&\frac{V(r_3,r_4,R)}{V(r_1+\varepsilon_0,r_2,R)} \\
=&\frac{\int_{R-r_4}^{R-r_3}\bigl(\frac{J_{R}(t)}{J_{R}(R)}\bigr)^{n-1}dt}
{\int_{R-r_2}^{R-r_1-\varepsilon_0}\bigl(\frac{J_{R}(t)}{J_{R}(R)}\bigr)^{n-1}dt} \nonumber\\
\geq& \frac{\int_{R-r_4}^{R-r_3}(z(R-t)-\frac{R-t}{R}z(R))^{n-1}dt} {\int_{R-r_2}^{R-r_1-\varepsilon_0}(z(R-t))^{n-1}dt},\nonumber
\end{align}
where we use Lemma \ref{lem3.1} in the last inequality.
For $t\in[R-r_4,R-r_3]$, we have $z(R-t)-\frac{R-t}{R}z(R)\geq z(r_{4})-\frac{r_{4}}{R}>0$, hence
\begin{align}
\frac{\mathrm{Vol}(\Omega(K)\cap b^{-1}([r_{3}-\varepsilon_0,r_{4}+\varepsilon_0]))}{\mathrm{Vol}(\Omega(K)\cap b^{-1}([r_{1}+2\varepsilon_0,r_{2}-\varepsilon_0]))}\geq \frac{(r_4-r_3)(z(r_{4})-\frac{r_{4}}{R})^{n-1}} {(r_{2}-r_1-\varepsilon_0)(z(r_{1}+\varepsilon_{0}))^{n-1}}.
\end{align}

For every fixed $k>k_{0}:=\sup_{x\in K}b^{-1}(x)$, take $r_i=k+i$, and let $\varepsilon_0\rightarrow0$ and $R\rightarrow+\infty$, we have
\begin{align}
  \frac{\mathrm{Vol}(\Omega(K)\cap b^{-1}([k+3,k+4]))}{\mathrm{Vol}(\Omega(K)\cap b^{-1}([k+1,k+2]))}
\geq\bigl(\frac{z(k+4)}{z(k+1)}\bigr)^{n-1}.
\end{align}
This finishes the proof of Claim 1.
\end{proof}

Now we go to prove (\ref{3.3}).
For every $k\in \mathbb{Z}^{+}$, denote by $B_{k}=\ln(\mathrm{Vol}(\Omega(K)\cap b^{-1}([k+3,k+4])))$.
According to (\ref{3.2}), for $k\geq k_{0}+4$, we have
\begin{align}
B_{k}&\geq B_{k-2}+(n-1)[\ln(z(k+4))-\ln(z(k+1))]\\
& = B_{k-2}-(n-1)\int_{k+1}^{k+4}u(s)ds.\nonumber
\end{align}
By induction, we have
\begin{align}
B_{k}\geq \min\{B_{k_{0}+1},B_{k_{0}}\}-3(n-1)\int_{0}^{\infty}u(s)ds.
\end{align}
Hence (\ref{3.3}) holds for $k\geq k_{0}$.
The proof is completed.
\end{proof}

Now we give the proof of Proposition \ref{b}.

\begin{proof}

We first consider the case that $\gamma(0)=p$.
By (a) of Lemma \ref{br} we know that $r(x)\geq b(x)$ for $x\in b^{-1}((0,+\infty))$, hence
\begin{align}
\mathrm{Ric}_g(x)\geq -f(r(x))\geq -f(b(x)).
\end{align}
Set $K=\bar{B}(q,1)$.
According to Proposition \ref{3.1},
there exist $c_{0}>0$ and $k_{0}>0$ such that
\begin{equation}\label{c0}
\mathrm{Vol}(\Omega(K)\cap b^{-1}([k,k+1]))>c_0
\end{equation}
for any $k>k_{0}$.

On the other hand, for any $y\in \Omega(K)\cap b^{-1}([k_{0},k_{0}+r])$, there exists $y'\in K$ such that $y=\gamma_{y'}\cap b^{-1}(t)$, for some $t\in[k_{0},k_{0}+r]$.
By Lemma \ref{sormani6}, it is easy to see that $d(y,y')\leq k_{0}+1+r$, and hence $d(y,p)\leq k_{0}+2+r$, and
\begin{equation}\label{bB}
\Omega(K)\cap b^{-1}([k_{0},k_{0}+r])\subset B(p,k_{0}+2+r).
\end{equation}

Combine (\ref{bB}) with (\ref{c0}), it holds that
\begin{equation}
\mathrm{Vol}(B(p,k_{0}+2+r))\geq \mathrm{Vol}(\Omega(K)\cap b^{-1}([k_{0},k_{0}+r]))\geq c_0(r-1).
\end{equation}

Then for sufficiently large $r$, we have
\begin{equation}
\mathrm{Vol}(B(p,r))>\frac{c_0}{2}r.
\end{equation}

If the given ray $\gamma$ does not emit from $p$, then we just need to consider a Busemann ray $\gamma_p$ associated with $\gamma$, and consider the Busemann function $b_{\gamma_p}$ associated with $\gamma_p$.
Then we use the fact that $b_{\gamma}^{-1}(k)=b_{\gamma_{p}}^{-1}(k)$ for every large $k>0$ to complete the proof.
\end{proof}

\section{Part (1) of Theorem \ref{thm-opti-vol}}\label{sec-4}

Recall the following basic lemma, which is a consequence of the theory of $\mu$-bubbles due to Gromov \cite{G18}.

\begin{lem}[Lemma 2.1 in \cite{CLS23}]\label{lem:mububble}
Let $(N^3, g)$ be a 3-manifold with boundary satisfying $\mathrm{Sc}_g \geq 1$. Then there are universal constants $L > 0$ and $c > 0$ such that if there is a $p \in N$ with $d_N(p, \partial N) > L/2$, then there is an open set $\Omega \subset B(\partial N,L/2) \cap N$ and a smooth surface $\Sigma^2$ such that $\partial \Omega = \Sigma \sqcup \partial N$ and each component of $\Sigma$ has diameter at most $c$.
\end{lem}

In the following we fix the positive constants $L$ and $c$ so that $L\gg c$ and they satisfies the conclusion in Lemma \ref{lem:mububble}.

\begin{lem}\label{lem:main}
Let $(M^3, g)$ be a complete noncompact $3-$dimensional manifold with $\mathrm{Sc}_g\geq2$ outside a compact set $K$. Suppose that $M$ has finitely many ends and its first Betti number $b_1(M)<\infty$. Fix a point $p\in M$, $r(x)=d(x,p)$.
Suppose there is a function $f:\mathbb{R}^{+}\rightarrow \mathbb{R}^{+}$, such that $\displaystyle\lim_{r\rightarrow\infty}f(r)=0$ and
\begin{equation}\label{4.1}
\mathrm{Ric}_g(x)\geq -f(r(x)).
\end{equation}
Then there is an $r_0(p,M,g) > 0$ and a universal constant $C>0$ so that the following holds.
For every end $E\subset M$ such that $E$ contains a unbounded connected component of $M\setminus B(p,r_{0})$, any $r'>r_0$, and $a_{1}, a_{2}\in [L, 2L]$, let $E_1$ and $E_2$ be the unbounded components of $E \setminus \overline{B}(p,r'+a_1)$ and $E \setminus \overline{B}(p,r'+a_1+a_2)$ respectively, then $\bar{E}_1\setminus E_2$ consists of one connected component, and
\begin{equation}
\mathrm{diam}(\bar{E}_1\setminus E_2) \leq C.
\end{equation}
In particular, we have $\mathrm{diam}(\partial E_1)\leq C$.
\end{lem}
Lemma \ref{lem:main} slightly generalizes the results in \cite{CLS23}, and its proof follows the same line as that of \cite{CLS23}.
We give a sketch of its proof as follows.

\begin{proof}[Sketch of the proof of Lemma \ref{lem:main}:]
Let $E_0$ be the unbounded component of $E\setminus\bar{B}(p,r')$.
Since $b_{1}(M)$ is finite, following Proposition 3.2 in \cite{CLS22}, we can prove that there is a $r_0(p,M,g)>0$ such that $\partial E_k (k=0,1,2)$ are connected when $r'>r_0$.
In particular, we conclude that $\bar{E}_1\setminus E_2$ consists of one connected component.
Furthermore, we assume $r_0(p,M,g)$ is large so that $\mathrm{Sc}_g\geq2$ on $M\setminus B(p,r_{0})$.

Let $\gamma$ be a ray emitting from $p$ so that $\gamma(t)\in E$ for every large $t$.
Then for all $r'>r_0$, $\gamma\cap \partial B(p,r')$ lies on the boundary of the unique unbounded component of $E\setminus B(p,r')$.
We apply Lemma \ref{lem:mububble} to $E_k$ for $k = 0,1$. 
Similar to Lemma 5.4 in \cite{CLS22}, we obtain a connected surface $\Sigma_k$ in $B(\partial E_k,L/2) \cap E_k$ with $\mathrm{diam}(\Sigma_k) \leq c$ so that $\Sigma_k$ separates $\partial E_k$ from $E$.
Take $t_k \in \mathbb{R}_+$ with $\gamma(t_k) \in\Sigma_k$ for $k=0,1$.
Since $d(\gamma(t_0),\gamma(t_1)) \leq a_1 + a_2 \leq 4L$, we have $d_g(x_0, x_1) \leq 4L + 2c$ for any $x_0 \in \Sigma_0$ and $x_1 \in \Sigma_1$.

As in \cite{CLS23}, denote by $D:=4L + 2c, b:=3c, A:=2\sqrt{b^2+D^2}$. Fix $R\in \mathbb{R}$ such that
\begin{equation}\label{eq.assumption.R}
R\geq A+4L,\quad \sqrt{b^2+(2R+D)^2}+1 <2D +2R.
\end{equation}
Then take $\delta \in (0,1)$ small enough such that
\begin{equation}\label{eq.assumption.delta}
\delta<\sqrt{b^2+D^2},\quad 14\delta + 6\sqrt{\delta (D+R)}<c,\quad 2c+22 \delta +6\sqrt{\delta (D+R)}<\frac{L}{2}.
\end{equation}
All the constants above are independent of $(M,g)$ and $r_0$.

Since (\ref{4.1}) holds for $\displaystyle\lim_{r\rightarrow\infty}f(r)=0$, by Cheeger-Colding's almost splitting theorem (\cite{CC96}), if $r_0$ is sufficiently large, there is a length space $(Y, d_Y)$ with
\begin{equation}\label{4.5}
d_{GH}(B(\gamma(t_1),D+R) \subset (M, d_g),\  B((y,0),D+R)\subset (Y \times \mathbb{R}, d_Y \otimes d_{\text{Eucl}})) < \delta.
\end{equation}

The above is the only one occasion where the almost nonnegative Ricci curvature condition is used.
Once we have obtained (\ref{4.5}), together with the geometric meanings of constants in (\ref{eq.assumption.R}) and (\ref{eq.assumption.delta}), the remaining argument in the proof of Claims 1, 2, 3 in \cite{CLS23} can be carried out.
In particular, we can prove that,
$\bar{E}_1\setminus E_2$ is contained in the tubular neighborhood $T_{A}(\gamma|_{|t-t_1|<D+R})$ (this is just Claim 3 in \cite{CLS23}).
Therefore, we have
\begin{equation}
\mathrm{diam}(\bar{E}_1\setminus E_2)\leq 2A+2D+2R.
\end{equation}
This finishes the proof.
\end{proof}

With Lemma \ref{lem:main}, we can finish the proof of (1) of Theorem \ref{thm-opti-vol}.
\begin{proof}[Proof of (1) in Theorem \ref{thm-opti-vol}]
Let $E^{(1)},\ldots ,E^{(k)}$ be the unbounded connected components of $M\setminus \bar{B}(p,r_{0}+2L)$, where $r_{0}$ is from Lemma \ref{lem:main}.

For each $s\in \mathbb{Z}^{+}$, let $E^{(i)}_{s}$ denote the unbounded connected components of $E^{(i)}\setminus \bar{B}(p,r_{0}+2L+sL)$.
Then for every $r$ with $r_{0}+(s+2)L \leq r<r_{0}+(s+3)L$, we have
\begin{equation}
B(p,r) \subset (M \setminus \bigcup_{i=1}^{k}E^{(i)}_{1})\bigcup (\bigcup_{i=1}^{k}\bigcup_{j=2}^{s+1}(E^{(i)}_{j-1} \setminus E^{(i)}_{j}))
\end{equation}
Because the Ricci curvature of the manifold is asymptotically nonnegative, there exists $K\geq0$ such that $\mathrm{Ric}_g\geq-(n-1)K$.
By Lemma \ref{lem:main} and the Bishop-Gromov volume comparison theorem, there exists $\tilde{C}>0$ such that
\begin{equation}
\text{Vol}({E}_{j-1}^{(i)}\setminus E_{j}^{(i)})\leq\tilde{C}.
\end{equation}
Thus,
\begin{align}
\text{Vol}(B(p,r)) &\leq \text{Vol}(M \setminus \bigcup_{i=1}^{k}E^{(i)}_{1})+k\tilde{C}s \\
 &\leq \text{Vol}(M \setminus \bigcup_{i=1}^{k}E^{(i)}_{1})+k\tilde{C}(\frac{1}{L}(r-r_0)-2).\nonumber
\end{align}
i.e. $(M,g)$ has at most linear volume growth.
\end{proof}


\begin{cor}\label{cor4.3776}
Under the assumption of (1) in Theorem \ref{thm-opti-vol}, let $E$ be an end of $M$, then there exist $C_{2}>0$ and $R_{0}>0$ such that
\begin{align}\label{4.10}
\mathrm{diam}(E\cap \partial B(p,R))\leq C_{2}
\end{align}
for every $R\geq R_{0}$.
\end{cor}

\begin{proof}
Without loss of generality, we assume $E$ is one of the unbounded connected components of $M\setminus \bar{B}(p,r_{0}+2L)$, where $r_{0}$ is from Lemma \ref{lem:main}.
For every $R\geq r_{0}+2L$, $E\cap \partial B(p,R)$ may consists of many connected components, which is denoted by $A_{0}, A_{1}, \ldots, A_{\alpha}$.
We assume $A_{0}$ is the boundary of the unique unbounded connected component of $E\setminus B(p, R)$.
Then by Lemma \ref{lem:main}, $\mathrm{diam}(A_{0})\leq C$.
For $i\in\{1, \ldots, \alpha\}$, let
$$\mu_{i}:=\sup\{s|A_{i} \text{ is contained in the unbounded component of }E\setminus B(p,s)\}.$$
Then $\mu_{i}\geq r_{0}+2L$.
We apply Lemma \ref{lem:main} to $r'=\mu_{i}-1-L$, $a_{1}=a_{2}=L$ (recall that $L\gg 1$) to conclude that
every two points in $A_{i}$ and $\partial B_{i}$ respectively has a distance at most $C$, and
$$\mathrm{diam}(\partial B_{i})\leq C.$$
where $B_{i}$ denotes the unbounded components of $E \setminus \overline{B}(p,\mu_{i}-1)$.
Then we have $\mu_{i}\geq R-C+1$.
Let $i_{0}\in\{1, \ldots, \alpha\}$ is one such that $\mu_{i_{0}}$ is minimal in $\{\mu_{1}, \ldots, \mu_{\alpha}\}$.
From the construction it is easy to see that every point $x\in A_{i}$ can be connected by a geodesic to a point $y\in \partial B_{i_{0}}$, and $d(x,y)\leq C$.
On the other hand, basing on Lemma \ref{lem:main}, by an induction argument we can prove that every two points in $A_{0}$ and $\partial B_{i_{0}}$ respectively has a distance at most $(\frac{C}{L}+1)C$.
Thus we conclude that for every $R\geq r_{0}+2L$, $E\cap \partial B(p,R)=\bigcup_{i=0}^{\alpha}A_{i}$ has diameter at most $(\frac{C}{L}+2)C$.
The proof is completed.
\end{proof}


\begin{cor}\label{cor4.3}
Under the assumption of (1) in Theorem \ref{thm-opti-vol}, let $\gamma(t)$ be a ray emanating from $p\in M$, and let $b(x)$ be the Busemann function associated to $\gamma(t)$, then there exists a large constant $R_{0}>0$ such that
\begin{equation}\label{diam b}
\mathrm{diam}(b^{-1}(R))\leq 3C_{2},\quad \forall R\geq R_{0}.
\end{equation}
\end{cor}

\begin{proof}
Let $E$ be the end of $M$ such that $\gamma(t)\in E$ for every large $t$.
Recall that we have proved (\ref{4.10}) for every sufficiently large $R$.
Then by (b) of Lemma \ref{br}, (\ref{diam b}) holds.
\end{proof}

\section{Proof of Proposition \ref{lem-non-col}}\label{sec-5}
\begin{proof}[Proof of Proposition \ref{lem-non-col}.]

It suffices to prove that (\ref{1.9}) holds for every $q\in E$, where $E$ is any fixed end of $M$.
Consider a ray $\gamma(t)$ emanating from $p$, such that $\gamma(t)\in E$ for every large $t$, and let $b(x)$ denote the Busemann function associated to $\gamma(t)$.

Now, suppose there exists a sequence of points $q_i \in E$ such that $\mathrm{Vol}(B(q_i,1))\rightarrow 0$ as $i\rightarrow\infty$.
Let $d(p,q_i)=R_i$, then $R_i\rightarrow\infty$.
Since
\begin{equation}
d(b^{-1}(R_k),b^{-1}(R_k+1))=1,
\end{equation}
for every sufficiently large $k$, we have
\begin{align}\label{5.2}
&\mathrm{diam}(b^{-1}([R_k,R_k+1]))\\
\leq& \mathrm{diam}(b^{-1}(R_k))+\mathrm{diam}(b^{-1}(R_{k+1}))+d(b^{-1}(R_k),b^{-1}(R_k+1))\nonumber\\
\leq &6C_{2} + 1,\nonumber
\end{align}
where we use (\ref{diam b}) in the last inequality.
Then combining (\ref{4.10}) (\ref{5.2}) with the facts that $q_k\in\partial B(p,R_k)$, $\partial B(p,R_k)\cap b^{-1}([R_k,R_k+1])\neq\varnothing $, we have
\begin{equation}
b^{-1}([R_k,R_k+1])\subset \bar{B}(q_k,7C_{2}+1).
\end{equation}

Because the Ricci curvature of the manifold is asymptotically nonnegative, there exists $K\geq0$ such that $\mathrm{Ric}_g\geq-(n-1)K$.
By the Bishop-Gromov volume comparison theorem, there exists a constant $\tilde{C}$ depending only on $C_{2}$ and $K$ such that
\begin{align}
  &\text{Vol}(B(q_k,1))\\
\geq& \tilde{C}\text{Vol}(B(q_k,7C_{2}+1))\nonumber\\
\geq &\tilde{C}\text{Vol}(b^{-1}([R_k,R_k+1])).\nonumber
\end{align}

On the other hand, by Proposition \ref{3.1},
\begin{equation}
\mathrm{Vol}(b^{-1}([R_k,R_k+1])) > C_1 > 0,
\end{equation}
holds for every sufficiently large $k$,
then
\begin{equation}
\mathrm{Vol}(B(q_k,1)) \geq \tilde{C}C_1,
\end{equation}
leading to a contradiction.
The proof is completed.
\end{proof}

\section{Part (2) of Theorem \ref{thm-opti-vol}}\label{sec-6}

Firstly, let's recall some results on Ricci flow.
The following theorem is from \cite{BCW17} and \cite{ST17}.
\begin{thm}\label{local-RF}
Given any $\alpha_{0}>0$, $v_0>0$, there exist positive constants $T=T(\alpha_{0},v_0)$, $C=C(v_0)$ such that the following holds.
Let $(M, g)$ be a 3-dimensional Riemannian manifold, $B_{g}(x_{0},s)\subset\subset M$ for some $s\geq4$.
Suppose
\begin{align}
\mathrm{Ric}_{g}\geq -\alpha_{0} \quad \text{on } B_{g}(x_{0},s),
\end{align}
\begin{align}
\mathrm{Vol}(B_{g}(x,1))\geq v_0 \quad \text{for all }x\in B_{g}(x_{0},s-1).
\end{align}
\end{thm}
Then there exists a Ricci flow $(g(t))_{t\in[0, T]}$ on $B_{g}(x_{0},s-2)$, with $g(0) = g$, so that
\begin{align}
\mathrm{Ric}_{g(t)} \geq -C\alpha_{0} \quad \text{on }B_{g}(x_{0},s-2),
\end{align}
\begin{align}
|\mathrm{Rm}_{g(t)}|\leq \frac{C}{t} \quad \text{on }B_{g}(x_{0},s-2),
\end{align}
for all $t\in(0, T]$.

Recall the shrinking balls lemma and the expanding balls lemma from \cite{ST16} as follows.
\begin{lem}[The shrinking balls lemma, Corollary 3.3 in \cite{ST16}]\label{shrink-lem}
Suppose $(M, g(t))_{t\in[0,T]}$ is a Ricci flow on a manifold $M$ of any dimension $n$.
Then there exists $\beta=\beta(n)\geq1$ such that the following is true.
Suppose $x_{0}\in M$ and that $B_{g(0)}(x_{0},r)\subset\subset M$ for some $r > 0$,
and $|Rm|_{g(t)}\leq\frac{c_{0}}{t}$ for some $c_{0}> 0$.
Then for any $0\leq s\leq t\leq T$, we have
\begin{align}
B_{g(0)}(x_{0}, r) \supset B_{g(s)}(x_{0}, r-\beta\sqrt{c_{0}s})\supset B_{g(t)}(x_{0}, r-\beta\sqrt{c_{0}t}).
\end{align}
\end{lem}

\begin{lem}[The expanding balls lemma, Lemma 3.1 in \cite{ST16}]\label{expand-lem}
Suppose $(M, g(t))$ is a Ricci flow for $t\in[-T, 0]$, $T > 0$, on a manifold $M$ of any dimension.
Suppose that $x_{0}\in M$ and that $B_{g(0)}(x_{0},R)\subset\subset M$ and $\mathrm{Ric}_{g(t)}\geq -K < 0$ on $B_{g(0)}(x_{0},R)\cap B_{g(t)}(x_{0},Re^{Kt})$ for each $t\in[-T, 0]$.
Then for all $t\in [-T, 0]$, we have
\begin{align}
B_{g(0)}(x_{0},R)\supset B_{g(t)}(x_{0},Re^{Kt}).
\end{align}
\end{lem}

The following theorem provides a local maximum principles along Ricci flow.

\begin{thm}[Theorem 1.1 in \cite{LT22}]\label{local-max-prin}
Let $(M^{n}, g(t))$, $t\in[0, T]$, be a smooth solution to the Ricci flow which is possibly incomplete.
Suppose
\begin{align}
\mathrm{Ric}_{g(t)}\leq \frac{\alpha}{t}
\end{align}
on $M\times(0, T]$ for some $\alpha> 0$.
Let $\varphi(x, t)$ be a continuous function on $M\times[0, T]$ which satisfies $\varphi(x, t)\leq \frac{\alpha}{t}$ on $M\times(0, T]$ and
\begin{align}
\bigl(\frac{\partial}{\partial t}-\Delta_{g(t)}\bigr)\varphi \bigl|_{(x_{0},t_{0})}\leq L(x_{0},t_{0})\varphi(x_{0},t_{0})
\end{align}
whenever $\varphi(x_{0},t_{0})> 0$ in the sense of barrier, for some continuous function $L(x, t)$ on $M\times[0,T]$ with $L(x, t)\leq \frac{\alpha}{t}$.
Suppose $p\in M$ such that $B_{g(0)}(p,2)\subset\subset M$ and $\varphi(x,0)\leq0$ on $B_{g(0)}(p,2)$.
Then for any $l>\alpha+1$, there exists $\hat{T}=\hat{T}(n,\alpha,l) > 0$ such that for $t\in [0, \min\{T,\hat{T}\}]$,
\begin{align}
\varphi(p,t)\leq t^{l}.
\end{align}
\end{thm}

\begin{cor}\label{Scal-bound}
Let $(M^{n}, g(t))$, $t\in[0, T]$, be a smooth solution to the Ricci flow which is possibly incomplete.
Suppose
\begin{align}
\mathrm{Ric}_{g(t)}\leq \frac{\alpha}{t}
\end{align}
on $M\times(0, T]$ for some $\alpha> 0$.
Suppose $p\in M$ such that $B_{g(0)}(p,2)\subset\subset M$ and $\mathrm{Sc}_{g(0)}(p)\geq 2$ on $B_{g(0)}(p,2)$.
Then there exists $\hat{T}=\hat{T}(n,\alpha) > 0$ such that for $t\in [0, \min\{T,\hat{T}\}]$,
\begin{align}\label{3.5}
\mathrm{Sc}_{g(t)}(p)\geq2-t^{n\alpha+3}.
\end{align}
\end{cor}

\begin{proof}
Take $\varphi(x, t)=2-\mathrm{Sc}_{g(t)}$, and $L(x, t)=\frac{4+2\mathrm{Sc}_{g(t)}}{n}$, then $\varphi(x, t)\leq \frac{n\alpha}{t}$, $L(x, t)\leq \frac{n\alpha+1}{t}$, and
\begin{align}
\bigl(\frac{\partial}{\partial t}-\Delta_{g(t)}\bigr)\varphi =-2|\mathrm{Ric}_{g(t)}|^{2}\leq -2\frac{\mathrm{Sc}^{2}_{g(t)}}{n} <L\varphi.
\end{align}
According to Theorem \ref{local-max-prin}, (\ref{3.5}) holds for suitable $\hat{T}=\hat{T}(n,\alpha)$.
\end{proof}


The following theorem, which will be used in our proof of Theorem \ref{thm-opti-vol}, slightly generalising Theorem 1.8 in \cite{ST17}.
Such kind of result may be known to experts, and its proof is a slight modification of that of \cite{ST17}.
We include the proof for completeness.

\begin{thm}\label{thm-RF-from-RLS}
Let $(M_{i}, g_{i}, p_{i})$ be a sequence of 3-dimensional Riemannian manifolds (possibly incomplete).
Suppose $U_{i}=B_{g_{i}}(p_i, s_{i})\subset\subset M_{i}$ for a sequence of $s_{i}\rightarrow\infty$, and  $\mathrm{Ric}_{g_{i}}\geq-\alpha_{i}$ (where $\alpha_{i}\in[0,1]$, and $\lim_{i\rightarrow\infty}\alpha_{i}=\alpha_{\infty}$) on $B_{g_{i}}(p_i, s_{i})$, and $\mathrm{Vol}_{g_{i}}B(p, 1) > v_{0} > 0$ for every $p\in B_{g_{i}}(p_i, s_{i}-1)$ and all $i\in\mathbb{N}$.
Then there exists a smooth manifold $M$, a point $p_{\infty}\in M$, a complete Ricci flow $g(t)$ on $M$ for $t \in(0, T]$, where $T> 0$ depends only on $\alpha_{\infty}$ and $v_{0}$, and a continuous distance metric $d_{0}$ on $M$ such that $d_{g(t)}\rightarrow d_{0}$ locally uniformly as $t\downarrow0$,
and up to passing to a subsequence we have that $(U_{i}, g_{i}, p_{i})$ converges in the pointed Gromov-Hausdorff sense to $(M, d_{0}, p_{\infty})$.
Furthermore, the Ricci flow satisfies
\begin{align}\label{6.13}
\mathrm{Ric}_{g(t)} \geq -C\alpha_{\infty} \quad \text{on }M\times(0,T],
\end{align}
\begin{align}\label{6.14}
|\mathrm{Rm}_{g(t)}|\leq \frac{C}{t} \quad \text{on }M\times(0,T],
\end{align}
\begin{align}\label{6.15}
\mathrm{inj}_{g(t)}\geq \rho_{0}\sqrt{t} \quad \text{for all }t\in(0,T],
\end{align}
\begin{align}\label{6.16}
d_{g(t_{1})}(x, y)-\beta\sqrt{C}(\sqrt{t_{2}}-\sqrt{t_{1}})\leq d_{g(t_{2})}(x, y)\leq e^{C\alpha_{\infty}(t_{2}-t_{1})} d_{g(t_{1})}(x, y)
\end{align}
for any $0\leq t_{1}\leq t_{2}\leq T$ and $x,y\in M$ (where $d_{g(0)}$ means $d_0$), where $C$ depends only on $v_0$, and $\beta> 0$ is a universal constant.
\end{thm}

\begin{proof}
Without loss of generality, we assume $\alpha_{i}\leq\alpha_{\infty}+1$ for every $i$. 
By Theorem \ref{local-RF}, there exists a Ricci flow $(g_{i}(t))_{t\in[0, T]}$ on $B_{g_{i}}(p_{i},s_{i}-2)$, with $g_{i}(0) = g$, so that
\begin{align}\label{4.133}
\mathrm{Ric}_{g_{i}(t)} \geq -C\alpha_{i} \quad \text{on }B_{g}(p_{i},s_{i}-2)\times(0,T],
\end{align}
\begin{align}\label{4.137}
|\mathrm{Rm}_{g_{i}(t)}|\leq \frac{C}{t} \quad \text{on } B_{g}(p_{i},s_{i}-2) \times(0,T],
\end{align}
for some $T=T(\alpha_{\infty},v_0)$, $C=C(v_0)$.

By lemmas \ref{shrink-lem} and \ref{expand-lem}, up to replacing $T$ by a smaller one, for any $0\leq s\leq t\leq T=T(\alpha_{\infty},v_0)$, we have
\begin{align}
B_{g_{i}(0)}(p_{i}, s_{i}-2) \supset B_{g_{i}(t)}(p_{i}, s_{i}-4),
\end{align}
\begin{align}\label{5.5}
B_{g_{i}(t)}(p_{i}, \frac{1}{2}s_{i}-2) \supset B_{g_{i}(s)}(p_{i}, \frac{1}{4}s_{i}-1).
\end{align}

Denote by $U_{i}'=B_{g_{i}(0)}(p_{i}, \frac{1}{4}s_{i}-2)$ for short.
For any $0\leq t\leq T$ and $x, y\in U_{i}'$, by (\ref{5.5}), the shortest curve with respect to $g_{i}(t)$ connecting $x$ and $y$ is contained in $B_{g_{i}(t)}(p_{i}, s_{i}-4)\subset B_{g_{i}(0)}(p_{i}, s_{i}-2)$.
Making use of (\ref{4.133}) and (\ref{4.137}), it is well-known that, there exists $\beta> 0$ such that for any $0\leq t_{1}\leq t_{2}\leq T$ and $x, y\in U_{i}'$ (see e.g. \cite{ST16}), we have
\begin{align}\label{5.60}
d_{g_{i}(t_{1})}(x, y)-\beta\sqrt{C}(\sqrt{t_{2}}-\sqrt{t_{1}})\leq d_{g_{i}(t_{2})}(x, y)\leq e^{C\alpha_{i}(t_{2}-t_{1})} d_{g_{i}(t_{1})}(x, y).
\end{align}

By (\ref{4.137}), according to Shi's derivative estimates, for any integer $k\geq1$,
\begin{align}\label{4.136}
|\nabla^{k} \mathrm{Rm}_{g_{i}(t)}|\leq \frac{C_{k}}{t^{1+\frac{k}{2}}} \quad \text{on } U_{i}',
\end{align}
where $C_{k}$ depends on $C$ and $k$.
By Lemma 4.1 in \cite{ST17}, there exists a positive constant $\rho_{0}$ depending only on $v_{0}$ so that
\begin{align}
\mathrm{inj}_{g_{i}(t)}\geq \rho_{0}\sqrt{t} \quad \text{on } U_{i}'.
\end{align}
By Hamilton's compactness theorem, up to choosing a subsequence, $(U_{i}', g_{i}(t), p_{i})_{t\in(0,T]}$ converges in the smooth Cheeger-Gromov sense to a Ricci flow solution $(M, g(t), p_{\infty})_{t\in(0,T]}$, where $M$ is a three-manifold.
By the shrinking balls lemma, it is easy to see that for any $t\in(0,T]$, $(M,g(t))$ is complete.
By the smooth convergence for $t\in(0,T]$, it is obvious that (\ref{6.13})-(\ref{6.15}) hold.
Note that (\ref{6.16}) holds for any $0< t_{1}\leq t_{2}\leq T$ and $x,y\in M$, which ensures the existence of  a metric $d_{0}$ on $M$ so that $d_{g(t)}$ converges to $d_0$ locally uniformly as $t\rightarrow0$, and for any $x,y\in M$ and $0< t\leq T$, we have
\begin{align}\label{6.19}
d_{0}(x, y)-\beta\sqrt{Ct}\leq d_{g(t)}(x, y)\leq e^{C\alpha_{\infty}t}d_{0}(x, y).
\end{align}
By (\ref{6.19}), (\ref{5.60}) and the Cheeger-Gromov convergence of $(U_{i}',g_{i}(t))$  to $(M,g(t))$ for every small $t>0$, by the triangle inequality it is easy to check that $(U_{i}',g_{i}(0))$ converges to $(M,d_{0})$ in Gromov-Hausdorff sense.
\end{proof}

\begin{cor}\label{cor-scal-bound}
If in Theorem \ref{thm-RF-from-RLS}, we assume in addition that $\mathrm{Sc}_{g_{i}}\geq 2$ on $U_{i}=B_{g_{i}}(p_i, s_{i})$, then on $(M,g(t))_{t\in(0,T]}$, we have
\begin{align}\label{3.52223512}
\mathrm{Sc}_{g(t)}\geq2.
\end{align}
\end{cor}
\begin{proof}
Since (\ref{4.137}) holds, according to Corollary \ref{Scal-bound}, for any $t\in(0,T]$ (here we can choose $T$ to a smaller one), we have
\begin{align}\label{3.551329}
\mathrm{Sc}_{g_{i}(t)}\geq2-t^{\tau}\quad \text{on } B_{g_{i}}(p_{i},s_{i}-4)
\end{align}
for some $\tau>0$ depending only on $v_{0}$.
Then by the pointed Cheeger-Gromov convergence of $(U_{i}',g_{i}(t))$  to $(M,g(t))$, we know
\begin{align}
\mathrm{Sc}_{g(t)}\geq2-t^{\tau}\quad \text{on } M\times (0,T].
\end{align}
For any $t_{0}\in (0,T]$, since the Ricci flow $(M,g(t))_{t\in[t_{0},T]}$ is complete and has uniformly bounded curvature, by the standard maximum principle, the lower bound of scalar curvature is preserved, i.e.
\begin{align}
\mathrm{Sc}_{g(t)}\geq2-t_{0}^{\tau}\quad \text{on } M\times [t_{0},T].
\end{align}
By the arbitrariness of $t_{0}$, we conclude that (\ref{3.52223512}) holds.
\end{proof}

\begin{proof}[Proof of (2) of Theorem \ref{thm-opti-vol}:]
By an argument by contradiction, we will prove that on every end $E$ with a ray $\gamma:[0,\infty)\rightarrow E$, it holds that
\begin{equation}\label{opti-vol-con}
\limsup_{r\rightarrow\infty}\frac{\mathrm{Vol}(E\cap B(p, r))}{r}\leq 4\pi.
\end{equation}
Let $b_{\gamma}$ be the Busemann function associated with $\gamma$.
We will prove that, for any $\epsilon>0$, there exists some $r_{0}>0$ so that for any $r>r_{0}$, it holds
\begin{align}\label{6.966}
\mathrm{Vol}(\{-1\leq b_{\gamma}-r\leq 1\})<8\pi+\epsilon.
\end{align}
If (\ref{6.966}) is proved, then for every $r>r_0$,
\begin{align}
&\mathrm{Vol}(E\cap B(p,r))\leq \mathrm{Vol}(E\cap b^{-1}([0,r])) \\
\leq& \mathrm{Vol}(E\cap b^{-1}([0,r_{0}]))+\mathrm{Vol}(E\cap b^{-1}([r_{0},r]))\nonumber\\
\leq& \mathrm{Vol}(E\cap b^{-1}([0,r_{0}]))+ (8\pi+\epsilon)(\frac{r-r_0}{2}+1), \nonumber
\end{align}
and hence $\limsup_{r\rightarrow\infty}\frac{\mathrm{Vol}(E\cap B(p, r))}{r}\leq 4\pi+\frac{\epsilon}{2}$.
Then (\ref{opti-vol-con}) holds by the arbitrariness of $\epsilon$.

We prove (\ref{6.966}) by contradiction.
Suppose on the contrary, there is some $\epsilon_{0}>0$ such that there is a sequence of $r_{i}\uparrow +\infty$ and
\begin{align}\label{6.966-con}
\mathrm{Vol}(\{-1\leq b_{\gamma}-r_{i}\leq 1\})\geq8\pi+\epsilon_{0}.
\end{align}
Denote by $p_{i}=\gamma(r_{i})$, $s_{i}=\frac{r_{i}}{2}$, $\epsilon_{i}=f(s_{i})$, $V_{i}=B_{g}(p_{i},s_{i})$.

Note that on $V_{i}$ the Ricci curvature has a lower bound $-\epsilon_{i}$ with $\epsilon_{i}\rightarrow 0$.
In addition, by Proposition \ref{lem-non-col}, there exists some $v_0>0$ such that $\mathrm{Vol}(B(p, 1)) > v_{0} > 0$ for all $p\in V_{i}$.
By Theorem \ref{thm-RF-from-RLS} and Corollary \ref{cor-scal-bound}, there exists a smooth manifold $X$, a point $p_{\infty}\in X$, a complete Ricci flow $g(t)$ on $X$ for $t \in(0, T]$, and a continuous distance metric $d_{0}$ on $X$ such that $d_{g(t)}\rightarrow d_{0}$ locally uniformly as $t\downarrow0$,
and up to passing to a subsequence we have that $(V_{i}, g, p_{i})$ converges in the pointed Gromov-Hausdorff sense to $(X, d_{0}, p_{\infty})$, and (\ref{3.52223512}) holds on $X\times(0,T]$.
Furthermore, for each $t\in(0, T]$, we have
\begin{align}\label{5.9}
\mathrm{Ric}_{g(t)} \geq 0, 
\end{align}
and for any $x,y\in X$ and $0< t\leq T$, we have
\begin{align}\label{5.8}
d_{0}(x, y)-\beta\sqrt{Ct}\leq d_{g(t)}(x, y)\leq d_{0}(x, y).
\end{align}

Note that $(X, d_{0}, \mathcal{H}^{3})$ is a non-collapsed $\mathrm{RCD}(0,3)$ space (we say $\mathrm{ncRCD}(0,3)$ for short) in the sense of \cite{GigPhi15}.
On the other hand, by the choice of $V_{i}$, it is obvious that $X$ contains a line, hence by the splitting theorem of Gigli (see \cite{Gig13}), there exists a $\mathrm{ncRCD}(0,2)$ space $(Y, d_{Y}, \mathcal{H}^{2})$ such that $(X, d_{0}, \mathcal{H}^{3})$ is isomorphic to $(\mathbb{R}\times Y, d_{\mathrm{Eucl}}\otimes d_{Y}, \mathcal{L}^{1}\otimes\mathcal{H}^{2})$.

Denote by $f_{i}=b_{\gamma}-r_{i}$, then each $f_{i}$ is $1$-Lipschitz with $f_{i}(p_{i})=0$.
According to Proposition 27.20 of \cite{Vi09}, up to a subsequence, $f_{i}$ converges locally uniformly to a $1$-Lipschitz function $P$ on $(X,d_0)$.
By Corollary \ref{cor4.3}, there exist positive constants $C_1$ and $R_{0}$ such that
\begin{align}\label{diambbound}
\mathrm{diam}(b_{\gamma}^{-1}(R))\leq C_{1}\quad \text{ for every } R\geq R_{0}.
\end{align}
Combining this fact with Lemma \ref{sormani6}, it is not hard to prove that
\begin{align}\label{6.88}
\mathrm{diam}(P^{-1}(R))\leq C_{1} \quad \text{ for every } R.
\end{align}
On the other hand, for any $\mu_{0}\in \mathbb{R}$, denote by $F_{i}:=f_{i}^{-1}([\mu_{0},\infty))$.
According to Lemma \ref{sormani6}, $\mu_{0}-f_{i}(x)=d(x,F_{i})$, for every $x\in V_{i}\setminus F_{i}$.
Since $\mu_{0}-f_{i}$ converges to $\mu_{0}-P$ uniformly, similar to Lemma 3.25 in \cite{H20}, one can prove that $\mu_{0}-P(x)=d(x,P^{-1}([\mu_{0},\infty))$ for every $x\in X$ satisfying $P(x)< \mu_{0}$.
By the arbitrariness of $\mu_{0}$, it is easy the see that there exists a line $\eta\subset X$ such that
\begin{align}\label{6.89}
P(\eta(t_{1}))-P(\eta(t_{2}))=t_{1}-t_{2}\quad \text{ for every }t_{1}, t_{2}\in \mathbb{R}.
\end{align}
Note that the projection of $\eta$ to the $Y$ factor must be a point.
For otherwise the projection on $Y$ is line, and by the splitting theorem of Gigli, $(X, d_{0})$ is isometric to $(\mathbb{R}^{2}\times Z, d_{\mathrm{Eucl}}\otimes d_{Z})$, contradicting to (\ref{6.88}) and (\ref{6.89}).
Thus the map $P$ is just the projection of $\mathbb{R}\times Y$ onto the $\mathbb{R}$-factor, and hence
\begin{align}\label{diamY}
\mathrm{diam}(Y)\leq C_{1}.
\end{align}

\textbf{Claim: }for every $t\in (0,T]$, $(X,g(t))$ contains a line.

\begin{proof}[Proof of the claim:] We will prove the claim holds for every sufficiently small $t>0$, then the Ricci flow solution ensures the validity of the claim for larger $t$.
Let $\eta:(-\infty, +\infty)$ be the line on $\mathbb{R}\times Y$ with $\eta(0)=p_{\infty}$.
Now we fix a small $t$.
For every large $R>0$, denote by $x_{\infty}=\eta(-R)$, $y_{\infty}=\eta(R)$.
By (\ref{5.8}), $|d_{g(t)}(x_{\infty}, y_{\infty})-2R|\leq C_{2}\sqrt{t}$, $|d_{g(t)}(p_{\infty}, x_{\infty})-R|\leq C_{2}\sqrt{t}$ and $|d_{g(t)}(p_{\infty}, y_{\infty})-R|\leq C_{2}\sqrt{t}$.
Let $z^{(R)}$ be the mid-point of a geodesic with respect to $g(t)$, $\eta^{(R)}:[-\frac{1}{2}d_{g(t)}(x_{\infty}, y_{\infty}), \frac{1}{2}d_{g(t)}(x_{\infty}, y_{\infty})]\rightarrow X$, which connects $x_{\infty}$ and $y_{\infty}$, then by (\ref{5.8}) again, $|d_{0}(x_{\infty}, z^{(R)})-R|\leq \frac{3}{2}C_{2}\sqrt{t}$ and $|d_{0}(y_{\infty}, z^{(R)})-R|\leq \frac{3}{2}C_{2}\sqrt{t}$.
Combining them with the product structure $\mathbb{R}\times Y$ and (\ref{diamY}), it is easy to see that there exists a $C_{3}>0$ depending only on $C_{1}$ and $C_{2}\sqrt{t}$ (we assume $R$ is sufficiently large, and it is easy to see that $C_{3}$ is independent of $R$) so that $d_{0}(p_{\infty}, z^{(R)})\leq C_{3}$.
Now we take a sequence of $R_{i}\rightarrow +\infty$, and take geodesics $\eta^{(R_{i})}$ with respect to $g(t)$, whose length $\approx 2R_{i}$, and with midpoint $z^{(R_{i})}$ contained in a compact set. After taking a converging subsequence, we obtain a line in $(X,g(t))$.
\end{proof}

By the above claim and (\ref{5.9}), $(M_{\infty},g_{\infty}(t))_{t\in (0,T]}$ is isometric to $(\mathbb{R},g_{\mathrm{Eucl}})\otimes (\tilde{S},\tilde{g}(t))_{t\in (0,T]}$, where $\tilde{g}(t)$, $t\in (0,T]$, is a Ricci flow solution on a smooth manifold $\tilde{S}$.
By (\ref{3.52223512}), it is easy to see that
\begin{align}
\mathrm{sec}_{\tilde{g}(t)}\geq1.
\end{align}
Hence $\tilde{S}$ is either diffeomorphic to $S^{2}$ or $RP^{2}$, and each $(\tilde{S},\tilde{g}(t))$ is a $\mathrm{ncRCD}(1,2)$ space.
Clearly, the pointed Gromov-Hausdorff convergence from $(X,g(t),p_{\infty})$ (where $t\downarrow 0$) to $(X, d_{0}, p)$ induces the Gromov-Hausdorff convergence from $(S^{2},\tilde{g}(t))$ to $(Y, d_{Y})$.
So the limit $(Y, d_{Y}, \mathcal{H}^{2})$ satisfies $\mathrm{ncRCD}(1,2)$ condition.
By the generalized relative volume comparison theorem for $\mathrm{CD}$-spaces (see \cite{St06II}), it is easy to see that
\begin{align}\label{6.14154}
\mathcal{H}^{2}(Y)\leq 4\pi.
\end{align}

Since $f_{i}=b_{\gamma}-r_{i}$ locally uniformly converges to the projection $P: \mathbb{R}\times Y\rightarrow \mathbb{R}$, and (\ref{diambbound}) (\ref{diamY}) hold, we have that (see e.g. Lemma 3.26 in \cite{H20} or Proposition 2.10 in \cite{Zxy23}),
for every large $i$,
\begin{align}
|\mathrm{Vol}(\{-1\leq b_{\gamma}-r_{i}\leq 1\})-\mathcal{H}^{3}(\{-1\leq P\leq 1\})|<\frac{\epsilon_{0}}{2}.
\end{align}
Combining it with (\ref{6.14154}), we have
\begin{align}
\mathrm{Vol}(\{-1\leq b_{\gamma}-r_{i}\leq 1\})<8\pi+\frac{\epsilon_{0}}{2},
\end{align}
contradicting to (\ref{6.966-con}).
This finishes the proof of Theorem \ref{thm-opti-vol}.
\end{proof}

\begin{rem}
The argument to derive (\ref{6.14154}) is motivated by \cite{ZZ23}, where the authors are concentrated on the optimal diameter upper bound of $Y$.
Some technical differences, such as the possible negative Ricci curvature and the non-global positive scalar curvature lower bound, require some adjustments in our proof.
In addition, our proof of the claim is different from \cite{ZZ23}.
In fact, the argument in \cite{ZZ23} to ensure the splitting of $(X,g(t))$ can be carried to our setting, but we provide such a different proof because we think that it may be used in other related problems.
\end{rem}

\end{document}